\documentclass[12pt]{article}

\usepackage{a4wide}
\usepackage{amssymb}
\usepackage{amsfonts}
\usepackage{amsmath}
\input xy
\xyoption{arrow} \xyoption{matrix}

\date{}

\newtheorem{proposition}{Proposition}[section]
\newtheorem{theorem}[proposition]{Theorem}

\newtheorem{corollary}[proposition]{Corollary}

\def\fdim{ {\rm fdim}\,}

\def\GK{{\rm  GK}\,}

\def\der{\partial }

\def\nFM0{{\nu }_{F,M_0}}
\def\nFN0{{\nu }_{F,N_0}}
\def\nGN0{{\nu }_{G,N_0}}

\def\N0{ {\bf N}_0 }

\def\t{\otimes}
\def\g{\gamma}

\def\ra{\rightarrow}

\def\Xpm{X^{\pm }}

\def\l1{{\lambda}_1}

\def\a{\alpha}
\def\a0{ {\alpha }_0}
\def\a1{ {\alpha }_1}

\def\l{\lambda}

\def\nFGM0{{\nu }_{F,G,M_0}}

\def\nFN0{{\nu}_{F,N_0}}

\def\sm{{\sigma}^m}

\def\sm1{{\sigma}^{-1}}

\def\smtp1{{\sigma}^{-t+1}}

\def\S1{S^{-1}}

\def\Xpm1{X^{\pm 1}_1}

\def\sPM1{{\sigma }^{\pm 1}}
\def\sMP1{{\sigma }^{\mp 1 }}

\def\di{{\rm d.ind}}

\def\CD{{\cal D}}

\def\Ytm1{Y^{t-1}}
\def\Yim1{Y^{i-1}}

\def\CF{{\cal F}}

\def\dim{{\rm dim }}

\def\SL2Z{ {\rm SL}_2({\bf Z}) }

\def\Gp1{ G^{1 , 1 } }
\def\P11{ P^{-1 , 1 } }
\def\Pp1{ P^{1 , 1 } }

\def\nCLsr{{}^\nu\kern-2pt {\cal L}^{\sigma , \rho  }}
\def\nP{{}^\nu \kern-2pt P}
\def\nL{{}^\nu\kern-2pt L}
\def\nLL{{}^\nu\kern-2pt \Lambda}
\def\nPsr{{}^\nu\kern-2pt P^{\sigma , \rho  }}
\def\nLsr{{}^\nu\kern-2pt L^{\sigma , \rho  }}
\def\nuCL{{}^\nu\kern-2pt  {\cal L}}
\def\nCLsr{{}^\nu\kern-2pt {\cal L}^{\sigma , \rho  }}
\def\nCL1m{{}^\nu\kern-2pt {\cal L}^{-1 , 1  }}
\def\x1nu{x^\frac{1}{\nu}}
\def\xm1nu{x^{-\frac{1}{\nu}}}

\def\ra{\rightarrow }

\def\nAM0{{\nu }_{{\cal A},M_0}}
\def\nAN0{{\nu }_{{\cal A},N_0}}

\def\di!{\frac{\der^i}{i!}}
\def\dik!{\frac{\der^k_i}{k!}}

\begin{document}

\author{V. V. \  Bavula and V. Hinchcliffe 
}

\title{Morita invariance of the filter dimension and of the inequality of Bernstein}

\maketitle
\begin{abstract}
It is proved that the {\em filter dimenion} is {\em Morita
invariant}. A direct consequence of this fact is the {\em Morita
invariance} of the {\em inequality of Bernstein}: {\em if an
algebra $A$ is Morita equivalent to the ring $\CD (X)$ of
differential operators on a smooth irreducible affine algebraic
variety $X$ of dimension $n\geq 1$ over a field of characteristic
zero then the Gelfand-Kirillov dimension} $ \GK (M)\geq n =
\frac{\GK (A)}{2}$
 {\em for all nonzero finitely generated $A$-modules $M$.} In
 fact, a more strong result is proved, namely, {\em a Morita
 invariance of the holonomic number for   finitely generated
 algebra}. As a direct consequence of this fact  an affirmative
  answer is given  to the question/conjecture posed by Ken Brown several
  years ago of whether an analogue
 of the inequality of Bernstein holds  for {\em the (simple) rational Cherednik
 algebras} $H_c$ {\em for integral} $c$: $\GK (M)\geq n =\frac{\GK
 (H_c)}{2}$ {\em for all nonzero finitely generated $H_c$-modules} $M$.

 {\em Mathematics subject classification 2000: 16P90, 13N10, 16S32, 16P90, 16D30, 16W70.}


\end{abstract}


\section{Introduction}
Throughout the paper, $K$ is a  field, a module $M$ over an
 algebra $A$  means a {\em left} module denoted ${}_AM$, and $\GK
 ({}_AM)$ is the {\em Gelfand-Kirillov dimension} of an $A$-module
 $M$.

 Intuitively, the filter dimension of an algebra $\fdim (A)$ or a
module $\fdim (M)$ measures how `close' {\em standard} filtrations
of the algebra or the module are. In particular, for a simple
algebra it also measures the growth of how `fast' one can prove
that the algebra is simple.

The filter dimension appears naturally when one wants to
generalize the inequality of  Bernstein for  the  ring $\CD (X)$
of differential operators (as above) to the class of {\em simple
finitely generated} algebras (recall that $\CD (X)$ is a {\em
simple finitely generated} algebra).

\begin{theorem}\label{IBerin}
{\rm (The inequality of Bernstein)} Let $X$ be a smooth
irreducible affine algebraic variety  of dimension $n= \dim
(X)\geq 1$ over a field of characteristic zero and $\CD (X)$ be
the ring of differential operators on $X$. Then
 $$ \GK (M)\geq n = \frac{\GK (\CD (X))}{2}$$
 for all nonzero finitely generated $\CD (X)$-modules $M$.
\end{theorem}

\begin{theorem}\label{InFFI}
\cite{Bavcafd} Let $A$ be a simple finitely generated algebra.
Then
$$
 \GK (M)\geq \frac{\GK (A)}{\fdim (A)+\max \{ \fdim (A),
1\} }
$$
 for all nonzero finitely generated $A$-modules $M$.
\end{theorem}

The following theorem is proved in Section \ref{FHMA}.

\begin{theorem}\label{25Apr06}
{\rm (Morita invariance of the filter dimension)} Morita
equivalent simple finitely generated algebras have the same filter
dimension.
\end{theorem}

As an immediate consequence one has the corollary.

\begin{corollary}\label{cIBerin}
{\rm (Morita invariance of the inequality of Bernstein)} Let $X$
and $\CD (X)$ be as in Theorem \ref{IBerin}. If the algebra $A$ is
Morita equivalent to the ring of differential operators $\CD (X)$
then
$$ \GK (M)\geq n = \frac{\GK (A)}{2}$$
  for all nonzero finitely generated $A$-modules $M$.
\end{corollary}

{\it Proof}. The algebra $\CD (X)$ is a simple finitely generated
algebra of Gelfand-Kirillov dimension $2n$, hence  the algebra $A$
is a simple finitely generated algebra of Gelfand-Kirillov
dimension $2n$ (since all the mentioned properties of algebra are
Morita invariant  and the algebra $A$ is Morita equivalent to $\CD
(X)$). The filter dimension $\fdim (\CD (X))=1$, \cite{Bavjafd}.
By Theorem \ref{25Apr06}, $\fdim (A)= \fdim (\CD (X))=1$, and by
Theorem \ref{InFFI}, $$ \GK (M)\geq \frac{2n}{1+1}=n $$
  for all nonzero finitely generated $A$-modules $M$.
$\Box $

\begin{corollary}\label{qa24May06}
Let $\CD (Q_c)$ be the ring of differential operators on
quasi-invariants (over $\mathbb{C}$) for integral values of $c$
(see \cite{BerEtGin03}, \cite{CV1}, and \cite{FV}). Then $\GK
(M)\geq n=\frac{\GK (\CD (Q_c))}{2}$ for all nonzero finitely
generated $\CD (Q_c)$-modules $M$.
\end{corollary}

{\it Proof}. The algebra $\CD (Q_c)$ is Morita equivalent to the
$n$'th Weyl algebra $A_n$, \cite{BerEtGin03}, and the result
follows from Corollary \ref{cIBerin}. $\Box $

For an algebra $A$,
$$ h_A:= \inf \{ \GK (M)\, | \, M\neq 0 \;\; {\rm is \; a \; finitely\;
generated\;}  A-{\rm module}\}$$ is called the {\em holonomic
number} of $A$. If $A$ is a simple finitely generated infinite
dimensional algebra then $h_A\geq 1$. For the ring $\CD (X)$  of
differential operators on a smooth irreducible affine algebraic
variety $X$ of dimension $n\geq 1$ over a field of characteristic
zero, $h_{\CD (X)}=n$. In particular, for the $n$'th Weyl algebra
$A_n$, $h_{A_n}=n$.
\begin{theorem}\label{25May06}
{\rm (Morita invariance of holonomic number)} Let $A$ and $B$ be
Morita equivalent finitely generated algebras. Then $h_A= h_B$.
\end{theorem}

{\it Proof}. This follows from Theorem \ref{M25May06} since
${}_AM\mapsto {}_BY\t_A M$ is the equivalence of the categories of
$A$-modules and $B$-modules. $\Box $

\begin{theorem}\label{M25May06}
Let $A$ and $B$ be Morita equivalent finitely generated algebras,
and  $\begin{pmatrix} A& X\\ Y & B\\
\end{pmatrix}$ be a Morita context that provides the Morita
equivalence ($XY=A$ and $YX=B$ where ${}_AX_B$ and ${}_BY_A$ are
bimodules (necessarily finitely generated  on {\em each side},
\cite{MR}, Sec. 3)). Then, for each finitely generated $A$-module
$M$, $\GK ({}_AM)= \GK ({}_BY\t_AM)$.
\end{theorem}

The next corollary gives an affirmative answer to the
question/conjecture posed by Ken Brown several years ago (in his
conference talk) of whether an analogue of the inequality of
Bernstein holds for (simple) rational Cherednik
 algebras $H_c$ for integral $c$.

\begin{corollary}\label{c25May06}
Let $H_c$ be the rational Cherednik algebra over $\mathbb{C}$
associated to a finite Coxeter group $W$ in an $n$-dimensional
vector space $h$ where $c$ is integral (see \cite{BerEtGin03} for
detail). Then $\GK (M)\geq n$ for all nonzero finitely generated
$H_c$-modules.
\end{corollary}

{\it Proof}. The algebra $H_c$ is Morita equivalent to the cross
product $A_n \# W$ where $A_n$ is the $n$'th Weyl algebra
\cite{BerEtGin03}, hence $h_{H_c}= h_{A_n \# W}$ (Theorem
\ref{25May06}). Since the action of the group $W$ on $A_n$ {\em
preserves} the natural (standard) filtration of the Weyl algebra
$A_n$  (by the total degree of the canonical generators $x_i$ and
$\der_j$), and the group $W$ is {\em finite}, it follows at once
from definition of the Gelfand-Kirillov dimension that $h_{A_n \#
W}=h_{A_n}=n$. Hence,
$$\GK (M)\geq  h_{H_c}=h_{A_n \#
W}=h_{A_n}=n$$ for all nonzero finitely generated $H_c$-modules
$M$.  $\Box $


\section{The filter dimension and the holonomic number are  Morita invariant}\label{FHMA}

{\bf The filter dimension}.  Let $\CF $ be the set of all
functions from the set of natural numbers $ \mathbb{N}=\{ 0, 1,
\ldots \}$ to itself. For each function $f\in \CF $, the
non-negative real number or $\infty $ defined as
$$ \g (f):=\inf \{  r\in \mathbb{R}\, | \, f(i)\leq i^r\; {\rm for
}\; i\gg 0\}$$ is called the  {\bf degree} of $f$. The function
$f$ has {\bf polynomial growth} if $\g (f)<\infty $. Let $f,g,
p\in \CF $, and $p(i)=p^*(i)$ for $i\gg 0$ where $p^*(t)\in
\mathbb{Q}[t]$ (a polynomial algebra with coefficients from the
field  of rational numbers). Then
\begin{eqnarray*}
\g (f+g)\leq \max \{ \g  (f), \g (g)\}, & & \g (fg)\leq \g (f)+ \g
(g),\\
\g (p)=\deg_t(p^*(t)), & & \g (pg)= \g (p)+ \g
(g).\\
\end{eqnarray*}
Let $A=K\langle a_1, \ldots , a_s\rangle $ be a finitely generated
 $K$-algebra. The finite dimensional filtration $F=\{ A_i\}_{i\geq 0} $
associated with the algebra generators $ a_1, \ldots , a_s$:
$$ A_0:=K\subseteq A_1:=K+\sum_{i=1}^sKa_i\subseteq \cdots \subseteq
A_i:=A_1^i\subseteq \cdots $$ is  called the {\em standard
filtration} for the algebra $A$. Let $M=AM_0$ be a finitely
generated $A$-module where $M_0$ is a finite dimensional
generating subspace. The finite dimensional filtration $\{
M_i:=A_iM_0\}$ is called the {\em standard filtration} for the
$A$-module $M$.

{\it Definition}. $\GK (A):=\g (i\mapsto \dim_K(A_i))$ and $ \GK
(M):=\g (i\mapsto \dim_K(M_i))$ are called the {\em
Gelfand-Kirillov} dimensions of the algebra $A$ and the $A$-module
$M$ respectively. For more detail on the Gelfand-Kirillov
dimension the reader is referred to the books \cite{KL},
\cite{MR}.

{\bf The return functions and the (left) filter dimension}.

 {\it Definition} \cite{Bavcafd}. The function $\nFM0 :\mathbb{N}\ra
\mathbb{N}\cup \{ \infty \}$,
$$ \nFM0 (i):=\min \{ j\in\mathbb{N}\cup \{ \infty \}: \;
A_jM_{i,gen}\supseteq M_0\;\; {\rm for \; all}\;\; M_{i,gen} \}$$
is called the {\em return function} of the $A$-module $M$
associated with the filtration $F=\{ A_i\}$ of the algebra $A$ and
the generating subspace $M_0$ of the $A$-module $M$ where
$M_{i,gen}$ runs through all generating subspaces for the
$A$-module $M$ such that $M_{i,gen}\subseteq M_i$.

{\it Definition}, \cite{Bavcafd}. $\fdim (M):= \g (\nu_{F, M_0})$
is called the {\em filter dimension} of the $A$-module $M$, and
$\fdim (A):= \fdim ({}_AA_A)$ is called the {\em filter dimension}
of the algebra $A$ (note that ${}_AA_A\simeq {}_{A\t A^{op}}A$
where $A^{op}$ is the {\em opposite} algebra to $A$).

It is easy to check that the filter dimension does {\em not}
depend  neither on the choice of the filtration $F$ nor $M_0$.

Suppose, in addition, that the finitely generated algebra $A$ is a
{\em simple} algebra.   The {\em return function} $\nu_F \in \CF $
and the {\em left return function} $\l_F\in \CF $ for the algebra
$A$ with respect to the standard filtration $F:= \{ A_i\}$ for the
algebra $A$ are defined by the rules:
\begin{eqnarray*}
\nu_F(i)&:=& \min \{ j\in \mathbb{N}\cup \{ \infty \} \, | \,\,
1\in A_jaA_j\;\, {\rm
for \; all}\; \,  0\neq a\in A_i\},\\
\l_F(i)&:=& \min \{ j\in \mathbb{N}\cup \{ \infty \} \, | \, 1\in
AaA_j\; \; {\rm for \; all}\; \; 0\neq a\in A_i\},
\end{eqnarray*}
where $A_jaA_j$ is the vector subspace of the algebra $A$ spanned
over the field $K$ by the elements $xay$ for all $x,y\in A_j$; and
$AaA_j$ is the left ideal of the algebra $A$ generated by the set
$aA_j$.  Under a mild restriction the maps $\nu_F(i)$ and
$\l_F(i)$ are finite (see \cite{FildimHazew}). It is easy to check
\cite{Bavcafd} that 
\begin{equation}\label{fdnug}
\fdim (A)= \g (\nu _F).
\end{equation}
More results on the filter dimension the reader can find in
\cite{FildimHazew}.

{\bf Morita equivalence and Morita context}. Let us recall some
facts we need on Morita equivalence, for  proofs the reader is
refereed to \cite{MR}: if algebras $A$ and $B$ are Morita
equivalent then $\GK (A) = \GK (B)$; $A$ is a simple algebra  iff
$B$ is a simple algebra; $A$ is a finitely generated algebra iff
$B$ is a finitely generated algebra.

Algebras $A$ and $B$ are {\em Morita equivalent} iff there is a
{\em Morita context} $\begin{pmatrix} A& X\\ Y & B\\
\end{pmatrix}$ with $XY=A$ and $YX=B$ where ${}_AX_B$ and
${}_BY_A$ are bimodules (necessarily finitely generated  on {\em
each side}, \cite{MR}, Sec. 3).

{\bf Proof of Theorem \ref{25Apr06}}.

Let $A$ and $B$ be simple finitely generated Morita equivalent
algebras.  The algebras $A$ and $B$ are Morita equivalent, hence,
by \cite{MR}, 3.5.4,
 there is a Morita context $\begin{pmatrix} A& X\\ Y & B\\ \end{pmatrix}$
with $XY=A$ and $YX=B$ where bimodules ${}_AX_B$ and ${}_BY_A$ are
finitely generated on {\em each side}, i.e. $X= AX_0= X_0B$ and
$Y=BY_0= Y_0A$ for some finite dimensional subspaces $X_0\subseteq
X$ and $Y_0\subseteq Y$. Let $\{ A_i\} $ and $\{ B_i\}$ be {\em
standard} filtrations on the algebras $A$ and $B$ respectively.
Enlarging, if necessary, the spaces $X_0$ and $Y_0$ one may assume
that there is a natural number $c$ such that
$$ A_1X_0\subseteq X_0B_c, \;\; X_0B_1\subseteq A_cX_0, \;\;
B_1Y_0\subseteq Y_0A_c, \;\; Y_0A_1\subseteq B_cY_0.$$ The
one-sided modules $X$ and $Y$ are equipped with following {\em
standard} filtrations:
$$ \{ X_i:= A_iX_0 \}, \;\; \{ X_i':= X_0B_i\}, \;\; \{ Y_i':=
B_iY_0 \} , \;\; \{ Y_i:= Y_0A_i\}.$$ Recall that $\fdim (A) = \g
( \nu )$ where
$$ \nu (i) :=\min \{ j \, | \, 1\in A_j a_iA_j\;\; {\rm for \;
all}\;\; 0\neq a_i\in A_i\},$$ where $A_j a_iA_j$ means the vector
subspace spanned by all the product $xa_ib$ for $x,y\in A_j$.
 Similarly, $\fdim (B) = \g (\nu')$ where $\nu' $ is defined as above
but for the algebra $B$ rather than $A$. For each $ i\geq 0$,
$$ A_iX_0\subseteq X_0B_{ci}, \;\; X_0B_i\subseteq A_{ci}X_0, \;\;
B_iY_0\subseteq Y_0A_{ci}, \;\; Y_0A_i\subseteq B_{ci}Y_0,$$ or,
equivalently,
$$X_i\subseteq X_{ci}', \;\; X_i'\subseteq X_{ci}, \;\;
Y_i'\subseteq Y_{ci}, \;\; Y_i\subseteq Y_{ci}'.$$ Here and
everywhere the product $V_1\cdots V_s$ of subspaces in an algebra
means the linear span of all the products $\{ v_1\cdots v_s\, | \,
v_1\in V_1, \ldots , v_s\in V_s\}$. $YX=B$ and $XY=A$ imply
\begin{equation}\label{3BCni}
1_B\in B_dY_0X_0B_d \;\; {\rm and}\;\; 1_A\in A_dX_0Y_0A_d
\end{equation}
for some natural number $d$. We may also assume that
$X_0Y_0\subseteq A_d$ and $Y_0X_0\subseteq B_d$. For each nonzero
element $a\in A_i$, $1_A\in A_{\nu (i)} aA_{\nu (i)}$, hence
\begin{equation}\label{1BCni}
1_B\in B_dY_01_AX_0B_d\subseteq  B_dY_0 \, A_{\nu (i)} aA_{\nu
(i)} \, X_0B_d\subseteq B_{ c\nu (i)+d} Y_0aX_0B_{ c\nu (i)+d}.
\end{equation}
For each $i\geq 0$ and each nonzero element $b\in B_i$, $X_0
B_dbB_dY_0\neq 0$ since, by (\ref{3BCni}),
$$ 0\neq b = 1_Bb1_B\in B_dY_0\cdot X_0B_dbB_dY_0\cdot X_0B_d.$$
Clearly, 
\begin{equation}\label{2BCni}
X_0B_dbB_dY_0\subseteq X_0B_{i+2d} Y_0\subseteq
A_{c(i+2d)}X_0Y_0\subseteq A_{c(i+2d) +d}.
\end{equation}
 Let $l(i):= c\nu
(ci+2cd +d)+d$. By (\ref{1BCni}) and (\ref{2BCni}), for each
nonzero element $b\in B_i$,
$$ 1_B\in B_{l(i)} Y_0\cdot X_0B_dbB_dY_0\cdot X_0B_{l(i)}\subseteq
B_{l(i)+2d} b B_{l(i)+2d}.$$ Therefore,
$$\nu'(i)\leq c\nu (c i+2cd +d) +3d, \;\; i\geq 0, $$
hence $\fdim (B)= \g (\nu' )\leq \g (\nu ) = \fdim (A)$. By
symmetry, we have the opposite inequality as well, and so $\fdim
(A) = \fdim (B)$. $\Box$

{\bf Proof of Theorem \ref{M25May06}}. We keep the notations of
the proof of Theorem \ref{25Apr06}. Enlarging, if necessary, the
spaces $X_0$ and $Y_0$ one may assume that $1_A\in X_0Y_0$ and
$1_B\in Y_0X_0$. Let $M= AM_0$ be a finitely generated $A$-module
($\dim (M_0)<\infty $) and $\{ M_i:= A_iM_0\}$ be the standard
filtration on $M$. The $B$-module $M':= Y\t_AM= BY_0M_0$ is a
finitely generated $B$-module with the finite dimensional
generating space $Y_0M_0$ which is the image of $Y_0\t_KM_0$ under
the map $\phi : Y\t_KM\ra Y\t_AM$. Similarly, for  vector
subspaces $V\subseteq Y$ and $U\subseteq M$, $VU$ stands for $\phi
(V\t_K U)$. Consider the standard filtration $\{ M_i':=
B_iY_0M_0\}$ on the $B$-module $M'$. Fix a natural number, say
$e$, such that $\dim (X_0), \dim (Y_0)\leq e$. First, we prove
that 
\begin{equation}\label{eMci}
\frac{\dim (M_i)}{e}\leq \dim (Y_0M_i) \leq \dim (M_{ci}'), \;\;
i\geq 0.
\end{equation}
$1_A\in X_0Y_0$ implies $M_i= 1_AM_i\subseteq X_0Y_0M_i$ and so
$\dim (M_i) \leq \dim (X_0)\dim (Y_0M_i)\leq e \, \dim (Y_0M_i)$,
this proves the first inequality. The second follows from the
following inclusions: $ Y_0M_i= Y_0A_iM_0\subseteq B_{ci} Y_0M_0=
M_{ci}'$.

By (\ref{eMci}), $\GK ({}_AM)\leq \GK ({}_BM')$. By symmetry, we
also have the inequality $\GK ({}_BM')\leq \GK ({}_AM'')$ where
$M'':= X\t_B Y\t_AM= AX_0Y_0M_0$, and $\{ M_i'':= A_iX_0Y_0M_0\}$
is the standard filtration on $M''$. Now, $1_A\in X_0Y_0$ implies
$M_i= A_i M_0= A_i1_AM_0\subseteq A_iX_0Y_0M_0= M_i''$. On the
other hand, $M_i''\subseteq X_0B_{ci} Y_0M_0\subseteq
X_0Y_0A_{c^2i} M_0= X_0Y_0M_{c^2i}$. Summarizing, we have
$$ \dim (M_i)\leq \dim (M_i'')\leq e^2\dim (M_{c^2i}), \;\; i\geq
0.$$ Hence, $\GK ({}_AM) = \GK ({}_AM'')$ which implies $\GK
({}_AM)= \GK ({}_BM')$ (since $\GK ({}_AM)\leq \GK ({}_BM')\leq
\GK ({}_AM'')= \GK ({}_BM)$),  as required. $\Box $

\begin{corollary}\label{c25Apr06}
Let $A$  be a simple finitely generated algebra and an algebra $B$
be  Morita equivalent to $A$. Then
$$ \GK (M)\geq h_A\geq \frac{\GK
(A)}{\fdim (A)+\max \{ \fdim (A), 1\} } $$
  for all nonzero finitely generated $B$-modules $M$.
\end{corollary}

{\it Proof}. The algebras $A$ and $B$ are Morita equivalent, hence
$\GK (A)= \GK (B)$, $B$ is a simple finitely generated algebra,
and $\fdim (A) = \fdim (B)$ (Theorem \ref{25Apr06}). By Theorems
\ref{InFFI} and \ref{25May06},
$$\GK (M)\geq h_B=h_A\geq \frac{\GK (B)}{\fdim (B)+\max \{ \fdim (B), 1\}
}=\frac{\GK (A)}{\fdim (A)+\max \{ \fdim (A), 1\} }$$
 for any nonzero finitely
generated $B$-module $M$.  $\Box $

\begin{corollary}\label{24May06}
Let $A$ be a simple finitely generated algebra, $e=e^2$ be a
nonzero idempotent of $A$, and $B:= eAe$. Then
$$\GK (M)\geq h_A\geq \frac{\GK (A)}{\fdim (A)+\max \{ \fdim (A), 1\} }$$
 for any nonzero finitely
generated $B$-module $M$.
\end{corollary}

{\it Proof}. The algebra $B$ is Morita equivalent to $A$,
\cite{MR} (the Morita context $\begin{pmatrix} A& Ae\\ eA & eAe\\
\end{pmatrix}$ provides a Morita equivalence since $eAAe= eAe$ and
$ AeeA=A$, since $A$ is simple and $e\neq 0$). Now, applying
Corollary \ref{c25Apr06} we have the result.  $\Box $

\begin{corollary}\label{c24May06}
Let $A$ be a simple finitely generated algebra, $G$ be a finite
group of algebra automorphisms of $A$ such that $|G|^{-1}\in A$,
$A^G$ be the fixed algebra. Then $$\GK (M)\geq h_A\geq \frac{\GK
(A)}{\fdim (A)+\max \{ \fdim (A), 1\} }$$
 for any nonzero finitely
generated $A^G$-module $M$.
\end{corollary}

{\it Proof}. The algebras $A$ and $A^G$ are Morita equivalent
\cite{MR}, 7.8.6. Now, use Corollary \ref{c25Apr06}.  $\Box $

$${\bf Acknowledgements}$$

The authors would like to thank Y. Berest, K. Brown, O. Chalykh,
I. Gordon and A.~Veselov for interesting discussions.

Department of Pure Mathematics

University of  Sheffield

Hicks Building

Sheffield S3~7RH

UK

email: v.bavula@sheffield.ac.uk

\end{document}